\documentclass[12pt,twoside]{article}
\usepackage[cp866]{inputenc}
\usepackage{amsfonts}
\usepackage{amssymb,amsmath}
\usepackage{amsthm}
\usepackage{euscript}
\usepackage[russian]{babel}
\usepackage{bm}

\makeatletter
\@addtoreset{footnote}{section}
\renewcommand{\section}{\@startsection{section}{1}%
{0pt}{3.5ex plus 1ex minus .2ex}%
{2.3ex plus.2ex}{\normalfont \Large \bfseries}}
\renewcommand{\subsection}{\@startsection{subsection}{1}%
{0pt}{3.25ex plus 1ex minus .2ex}%
{1.5ex plus .2ex}{\normalfont \large \bfseries}} \makeatother
\makeatletter
\renewcommand{\section}{\@startsection{section}{1}%
{20pt}{3.5ex plus 1ex minus .2ex}%
{2.3ex plus.2ex}{\normalfont \Large \bfseries}}
\makeatother
\begin{document}
\setcounter{page}{1}\numberwithin{equation}{section}

\newpage
\begin{center}
{\Large \bf
 Обзор результатов и открытых проблем по математическим моделям движения
вязкоупругих сред типа Джеффриса}
\\[15pt]
\large \bf Д.А. Воротников\footnote{CMUC,
Apartado 3008, 3001 - 454 Coimbra, Portugal,
mitvorot@mat.uc.pt}, В.Г. Звягин\footnote{Faculty of Mathematics, Voronezh State
University, Universitetskaya Pl. 1, 394006, Voronezh,
Russia, zvg@math.vsu.ru}
\end{center}

\newcommand {\R} {\mathbb{R}}
\newcommand {\E} {\mathbf{E}}
\def\be{\begin{equation}}
\def\ee{\end{equation}}
\def\fr#1#2{\frac{\partial #1}{\partial #2}}
\large

\begin{abstract} The Jeffreys model (also associated with the names of Lethersich and Oldroyd) is one of the crucial conceptions in the theory of viscoelastic fluids.  The models of Jeffreys type describe behaviour of  bitumens, blood, polymers and their solutions, dough, the earth's crust, concrete, lubricants etc. Study of BVPs corresponding to their statics and dynamics meets a lot of mathematical difficulties, which  turn out to be much harder than the ones that are related to the celebrated Navier-Stokes system. In this work, we make an attempt to review the recent results and main unsolved problems for equations of motion for the mediums of Jeffreys' type.

\end{abstract}

\section*{Введение}

Модель Джеффриса (называемая также именами Летерзиха или Олдройда) является центральной в теории вязкоупругих жидкостей. Она и связанный с ней класс моделей описывает поведение таких сред, как битумы, кровь, полимеры и их растворы, тесто, земная кора, бетон, смазки и др. Исследование краевых, начальных и начально-краевых задач, описывающих статику и динамику этих моделей сопряжено с больщими математическими трудностями, превышающими таковые для известной системы Навье-Стокса, которая также порождает ряд нерешенных математических проблем, начиная с проблемы глобального существования гладкого решения в случае размерности области более двух.

Данная работа посвящена обзору современных результатов
и основным нерешенным проблемам по уравнениям движения сред типа Джеффриса. Отметим, что имеется ряд недавних обзоров \cite{zv,zvobs,carar,galdi,gs2,YKwon,renar,seq,gran}, касающихся этой темы, но посвященных более общим или смежным вопросам, а потому не раскрывающих эту тематику должным образом.

Опишем вначале рассматриваемую ситуацию.

Итак, пусть $\Omega \subseteq {\mathbb{R}}^3$ - ограниченный сосуд
с границей $\Gamma = \partial \Omega$, целиком заполненный
некоторой средой, которую в дальнейшем будем называть жидкостью.
Жидкость будем представлять как совокупность материальных частиц,
заполняющих сосуд $\Omega$, причём эти частицы будем считать
настолько малыми, что их можно отождествлять с точками объёма
$\Omega$.

Под движением жидкости мы будем понимать движение материальных
точек объёма $\Omega$. Таким образом, описать путь, который
проделывает каждая точка объёма $\Omega$ за время $t_0\leq t \leq
T$, это и означает описать движение жидкости за это время.

Пусть в трёхмерном пространстве зафиксирована ортогональная
система координат и $e_1,e_2,e_3$ - векторы соответствующего
базиса. Сосуд с жидкостью $\Omega$ будем рассматривать как область
в трёхмерном пространстве~, а положение движущейся точки объёма
$\Omega$ ( или, что то же, частицы жидкости ) можно описывать с
помощью вектор-функции

$$x(t)=x_1 (t) e_1 + x_2 (t) e_2 + x_3 (t) e_3.$$

Ясно, что если в начальный момент  $t_0$  частица жидкости
занимала положение $x_0$, а её движение описывается с помощью
$x(t)$, то $x(t_0 )= x_0$.

Каждой частице объёма $\Omega$ соответствует  своя вектор-функция
$x(t)$, описывающая её движение. Движение жидкости будет описано,
если будут найдены все эти вектор-функции $x(t)$.

Зафиксируем момент времени $t$. В этот момент времени частица
жидкости, двигающаяся по закону $x(t)$, имеет скорость

$$\dot x(t)=\dot{x_1}(t) e_1 + \dot{x_2}(t) e_2 + \dot{x_3}(t) e_3.$$

Обозначим через $v(t,x)$  - скорость частицы жидкости,
находящейся в момент времени $t$ в точке $x$. Тогда

$$\dot x(t)=v(t,x(t)).$$

Отсюда следует, что если известны скорость движущейся жидкости в
каждой точке $x \in \Omega$ в каждый момент времени t, т.е.
известна вектор-функция $v (t,x)$, определённая для всех
$x\subset \Omega$ и $t\in [t_0,T]$, то, для того чтобы найти
вектор-функцию $x(t)$, описывающую движение частицы жидкости,
занимающей в начальный момент  $t_0$, положение $x_0$, надо решить
следующую задачу Коши для векторного дифференциального уравнения:

$$\frac{dx(t)}{dt} =v (t,x(t)),~~{}{}{} {}
   t\geq t_0 , ~x\in \Omega,$$
$$x(t_0)=x_0.$$

Таким образом, для того, чтобы описать движение жидкости,
достаточно знать распределение скоростей жидкости в каждой точке
$x\in \Omega$ и в каждый момент времени $t\in [t_0,T]$, или, что
то же, знать вектор-функцию $v(t,x)$.

Для определения вектор функции $v(t,x),$ исходя либо из
принципа Даламбера, либо из второго закона Ньютона, несложно
выводится (см. \cite{dya}) следующее уравнение, называемое общим уравнением
движения среды и записываемое ниже в векторной форме:

\begin{equation} \rho \left( \frac{\partial v}{\partial t} +
\sum\limits_{i=1}^{3} v_i \frac{\partial v}{\partial x_i}\right) -
Div T_H = \rho f.
\end{equation}

Здесь $f(t,x)$ - плотность внешних сил, действующих на  частицы
среды (например, плотность поля гравитационных сил), $ \rho (t,x)$
- плотность жидкости, $T_H$ -- тензор напряжений, а $Div T_H(t,x)$
-- дивергенция тензора $T_H (t,x)$, т.е. вектор

$$\begin{pmatrix}
\sum\limits_{j=1}^3\frac{\partial T_{H1j}(t,x)}{\partial x_j}\\
\sum\limits_{j=1}^3\frac{\partial T_{H2j}(t,x)}{\partial x_j}\\
\sum\limits_{j=1}^3\frac{\partial T_{H3j}(t,x)}{\partial
x_j}\end{pmatrix} =\sum\limits_{j=1}^3
\begin{pmatrix} \frac{\partial T_{H1j}(t,x)}{\partial x_j}\\
{}\\
\frac{\partial T_{H2j}(t,x)}{\partial x_j}\\
{}\\
\frac{\partial T_{H3j}(t,x)}{\partial x_j}\end{pmatrix}.$$

Для несжимаемой жидкости c постоянной плотностью\footnote{Ниже мы, как правило, будем касаться именно таких жидкостей, если не будет оговорено обратное.} имеется несколько больше информации.
Во-первых, для нее

\begin{equation} div
\,v(t,x)=\sum\limits_{j=1}^3 \frac{\partial v_j(t,x)}{\partial
x_j}=0.
\end{equation}

 Во-вторых, след тензора напряжений $T_H$
считается полностью независимым от характеристик деформации, и в
этом случае удобно ввести еще одну неизвестную скалярную функцию
$p(t,x)=-\frac{1}{3}Tr T_H$, называемую давлением, и тензор

\begin{equation} \sigma =T_H - \frac{1}{3}Tr T_H I,
\end{equation}
называемый иногда тензором касательных напряжений\footnote{Хотя это не совсем корректный термин}. Он характеризует силы
внутреннего трения в жидкости (здесь $I$ - единичный тензор).

В-третьих, будем считать плотность постоянной и равной единице.

Заметим, что $Div (pI)=\nabla p$. Поэтому (0.1) с учетом сделанных
замечаний и (0.3) принимает вид

\begin{equation} \frac{\partial v}{\partial t} +
\sum\limits_{i=1}^3 v_i \frac{\partial v}{\partial x_i} - Div
\sigma + \nabla p = f. \end{equation}

Это уравнение также называют общим уравнением движения несжимаемой
среды с постоянной плотностью.

Самым распространенным и физически обусловленным граничным условием для этих уравнений является условие прилипания частиц жидкости к границе области (или, более общо, что скорость стремится к нулю при приближении к границе):

\begin{equation} v|_{\partial \Omega} = 0. \end{equation}

Система (0.2),(0.4) описывает течение всех видов несжимаемых сред
с постоянной плотностью, но при этом она содержит девиатор тензора
напряжений, который явно не выражен через неизвестные $v$ и
$p$ этой системы. Чтобы выразить девиатор тензора напряжений
$\sigma$ через неизвестные системы (0.2),(0.4), как правило,
используют соотношения между девиатором тензора напряжений
$\sigma$ и тензором скоростей деформации $\mathcal{E} =
(\varepsilon_{ij}),$ где

 $$\varepsilon_{ij} = \varepsilon_{ij}(v) =
 \frac{1}{2}  \left(\frac{\partial v_i}{\partial x_j} + \frac{\partial v_j}{\partial
 x_i} \right).$$

 Отметим, что устанавливая соотношение между девиатором тензора
 напряжений и тензором скоростей деформации, тем самым
 устанавливается тип среды. Такие соотношения называются
 определяющими или реологическими соотношениями. Реологические
 соотношения относятся к разряду гипотез, которые устанавливаются,
 как правило, на основе экспериментов или методом механических
 моделей  (по поводу последнего метода более подробно см. ниже)
 и должны подтверждаться для конкретных жидкостей в результате
 опытных исследований.

 В течении последних полутора столетий основным объектом
 исследований математиков  в области гидродинамики является модель
 ньютоновской жидкости. Её реологическое соотношение (в "несжимаемой" ситуации) имеет вид:$$\sigma = 2\nu
\mathcal{E},~~\nu > 0,$$ где~ $\nu$~ -~ кинематический коэффициент
вязкости. Если мы подставим это соотношение в (0.4), то, учитывая
(0.2), получим хорошо известную систему уравнений Навье-Стокса:

 $$\frac{\partial v}{\partial t} + \sum\limits_{i=1}^3 v_i
\frac{\partial v}{\partial x_i} - \nu\Delta v + \nabla p = f.$$

Она описывает течение при умеренных скоростях большого числа
встречающихся на практике вязких несжимаемых жидкостей. Но уже в
середине $XIX$ века стали известны такие жидкости, которые не
подчиняются ньютоновскому определяющему соотношению. Таковыми
являются, например, жидкости, в которых после прекращения движения
напряжения не обращаются мгновенно в нуль, а спадают по некоторому
закону, то есть имеет место релаксация напряжений; а также
жидкости, в которых после снятия напряжений движение не
прекращается мгновенно, а затухает по некоторому закону, то есть
имеет место запаздывание деформаций; и жидкости, в которых имеют
место оба этих эффекта. Различные модели таких сред, учитывающих
предысторию движения, были предложены в $XIX$ веке Дж. Максвеллом,
В.Кельвином, В.Фойгтом, а в $XX$ веке Джеффрисом \cite{jeff}, Летерзихом \cite{rei}, Олдройдом \cite{old}, Ларсоном \cite{bookour}
и другими.

Первые два пункта данной работы связаны с моделированием вязкоупругих сред, в третьем дан общий обзор математических исследований по моделям типа Джеффриса, а в последующих ряд результатов, указанных в третьем пункте, изложен более детально. В последнем пункте обсуждаются ключевые нерешенные проблемы.

\section{Принцип объективности}

Это один из основных принципов рациональной механики, который
выражает тот факт, что свойства материала не зависят от выбора
наблюдателя. На основе этого принципа получили обоснование как ряд
известных ранее моделей, так возникли и их обобщения. Опишем
кратко этот принцип и факты, которые из него вытекают.

Наблюдателя в рациональной механике отождествляют с системой
отсчета, т.е. неким правилом, сопоставляющим каждой точке
пространства элемент $x$ пространства $R^3$, а каждому моменту
времени элемент $t$ числовой оси. При изменении наблюдателя
считается, что сохраняются метрики в $R^3$ и на числовой оси, и
сохраняется направление времени. Тогда самое общее изменение
координат каждой точки имеет вид
\begin{equation} t^*=t+a,\end{equation}
\begin{equation}x^*=x^*_0(t) +Q(t)(x-x_0), \end{equation}
где $a$ - некоторое значение времени, $x_0$ - некоторая точка в
пространстве, $x^*_0(t)$ - некоторая функция времени со значениями
в точках пространства, $Q$ - функция времени, значениями которой
являются ортогональные тензоры.

Замена наблюдателя индуцирует некоторое преобразование векторов и
тензоров. Принцип объективности утверждает, что формулы,
выражающие физические свойства тела и содержащие время $t$, точку
$x$ и их различные функции, не должны меняться при преобразованиях
(1.1)-(1.2).

Пусть имеется вектор, являющийся геометрическим, т.е.
представляющим из себя направленный отрезок, существующий
независимо от наблюдателя. Пусть в первой системе отсчета он имеет
вид $w=\overrightarrow{x_1x_2}$. Тогда в новой системе отсчета
$w^* =
\overrightarrow{x^*_1x^*_2}=\overrightarrow{x_0^*(t)+Q(t)(x_1-x_0),\,x_0^*(t)+Q(t)(x_2-x_0)}=Q(t)
\overrightarrow{x_1x_2}=Q(t)w$

Если же  $T$ - тензор, переводящий геометрические векторы в
геометрические, то имеем в первой системе отсчета: $$ w_1=Tw_2,$$
где $w_1$ и $w_2$ - два пробных геометрических вектора. Так как
для ортогонального тензора $Q(t)$ имеет место равенство $Q(t)^T
Q(t) =I$, то получаем
$$w^*_1=Q(t)w_1=Q(t)TQ(t)^TQ(t)w_2=T^*w_2^*, $$ где
$T^*=Q(t)TQ(t)^T.$

Исходя из выше написанного, будем называть векторнозначную функцию
$w$ времени и точки не зависящей от наблюдателя, если ее
координаты при изменении системы отсчета (1.1)-(1.2) преобразуются
следующим образом: $$ w^*(t^*,x^*)=Q(t)w(t,x),$$ а тензорнозначную
функцию $T$ не зависящей от наблюдателя, если справедливо
равенство:
$$T^*(t^*,x^*)=Q(t)T(t,x)Q(t)^T.$$ Из механики известно,
что тензор напряжений $T_H$ является таковым.

Скалярная функция $A$ времени и точки называется не зависящей от
наблюдателя, если при изменении системы отсчета она остаётся
неизменной:
$$A^*(t^*,x^*)=A(t,x).$$

Примером такой функции является плотность $\rho(t,x)$.

Далее, наряду с одним из важнейших тензоров гидродинамики,
тензором скоростей деформации
 $\mathcal{E}=\mathcal{E}(v)=(\varepsilon_{ij})=(\varepsilon_{ij}(v)),$
 в гидродинамике используется и тензор, называемый тензором завихренности $W(t,x)=(w_{ij}(v))$, где
$$ w_{ij} (v) =\frac{1}{2}  \left (\frac{\partial v_i}{\partial x_j}
- \frac{\partial v_j}{\partial x_i} \right).$$

 Таким образом, тензор скоростей деформации $\mathcal{E}$ является
 симметрической частью матрицы Якоби $\left(\frac{\partial v_i}{\partial
 x_j}\right)$, а тензор завихренности $W$ кососимметрической частью этой
 матрицы.
 \par Для проверки реологических соотношений на объективность используется
следующее утверждение о преобразовании при изменении наблюдателя
тензоров скоростей деформации и завихренности.

{\bf Теорема Зарембы-Зоравского.} При изменении системы отсчета
(1.1)-(1.2) имеют место следующие преобразования:
\begin{equation}
\mathcal{E}^*(t^*,x^*)=Q(t)\mathcal{E}(t,x)Q(t)^T,\end{equation}

$$W^*(t^*,x^*)=Q(t)W(t,x)Q(t)^T+Q'(t)Q(t)^T,$$
т.е. тензор скоростей деформации не зависит от наблюдателя, а
тензор завихренности от него зависит.

Как мы уже отмечали, тензор напряжений $T_H$ не зависит от
наблюдателя. Из этого следует, что тензор касательных напряжений
$\sigma$ также не зависит от наблюдателя. В самом деле,
$$\sigma ^*= T^*_H-\frac{1}{3}TrT^*_HI=$$
$$= Q(t)T_HQ(t)^T-\frac{1}{3}Tr(Q(t)T_HQ(t)^T)I=Q(t)T_HQ(t)^T-$$ $$-\frac{1}{3}\sum\limits_{i,j,k=1}^3
q_{ij}T_{Hjk}q_{ik}I
=Q(t)T_HQ(t)^T-\frac{1}{3}\sum\limits_{j,k=1}^3[Q(t)^TQ(t)]_{jk}T_{Hjk}I=$$
$$ = Q(t)T_H Q(t)^T-\frac{1}{3}\sum\limits_{j,k=1}^3\delta_{jk}T_{Hjk}I= $$
$$ =Q(t)T_HQ(t)^T-\frac{1}{3}TrT_HI=Q(t)\sigma Q(t)^T.$$
 Отсюда, в частности, имеем представление
\begin{equation}\sigma =Q(t)^T\sigma ^*Q(t).\end{equation} А из (1.3) следует, что
\begin{equation}\mathcal{E}=Q(t)^T\mathcal{E}^*Q(t).
\end{equation}

\section{Базовые модели движения вязкоупругих жидкостей} Для определения
реологических соотношений, описывающих движение вязкоупругих
жидкостей, часто применяется метод механических моделей. Кратко
опишем сущность этого метода.

Основные свойства для реологии --- упругость, вязкость и
пластичность. Для представления упругости используется спиральная
пружина. Для нее имеет место закон Гука: удлинение пружины прямо
пропорционально приложенной к ее концам силе. Эту модель будем
обозначать буквой $H$. Для представления вязкости используется
модель в виде пробирки, заполненной вязким маслом, в которой
свободно перемещается поршень. Скорость поршня (относительно
масла) прямо пропорциональна приложенной силе. Эта модель
обозначается $N$. Модель для пластичности нам не потребуется.

Для построения моделей тел со сложными реологическими свойствами
можно соединить эти элементы параллельно (обозначается $|$ ) или
последовательно (обозначается $-$). При параллельном  соединении
нагрузки, воспринимаемые каждым элементом, складываются, а
скорости удлинения каждого элемента одинаковы. При
последовательном соединении складываются скорости удлинения
элементов, и каждый из них подвергается одинаковой нагрузке.

Для реологических соотношений нужны не силы и скорости, а
напряжение и скорости деформации. Здесь нужно вспомнить, что нашей
целью является моделирование реальных вязкоупругих жидкостей.
Грубо говоря, эти жидкости предполагаются состоящими из
микрокомплексов из маленьких пружинок и пробирок с поршнями.
Поэтому мы должны понимать напряжение и скорость деформации так,
чтобы это согласовалось с соответствующими определениями этих
понятий для реальной жидкости. Исходя из этого, определим
напряжение $\sigma $ (в пружине или в пробирке с поршнем) как
отношение силы сопротивления (которая равна по модулю приложенной
силе) к площади поперечного сечения (пружины или поршня), а
скорость деформации $\mathcal{E} $ как половину отношения скорости
(поршня в масле или изменения длины пружины) к характерной
(средней) продольной длине пробирки с поршнем или пружины.

Теперь необходимо записать соотношения между напряжением и
скоростью деформации для пружины и пробирки с поршнем. Для
пробирки с поршнем это совсем просто: так как сила пропорциональна
скорости, то напряжение пропорционально скорости деформации

\begin{equation} \sigma _N = 2\eta \mathcal{E} _N .\end{equation}

Отметим, что здесь $\eta $ имеет физический смысл вязкости.

Условимся для механических моделей обозначать производную по
времени точкой. Продифференцировав по времени закон Гука для
пружины, получим, что производная от силы пропорциональна скорости
изменения длины пружины. Поэтому производная от напряжения
пропорциональна скорости деформации \begin{equation}\dot{\sigma
}_H=2\mu \mathcal{E} _H. \end{equation} Коэффициенты
пропорциональности $\mu$ и $\eta $ обычно положительны.

\par Наиболее простой конструкцией модели вязкоупругой жидкости
является тело Максвелла (с символьной записью $M=H-N$):
последовательно соединенные пружина и поршень. Выведем
реологическое соотношение для тела Максвелла. Для этого вспомним,
что при последовательном соединении напряжение постоянно
\begin{equation}\sigma _M=\sigma _N=\sigma _H,
\end{equation} а скорости деформации складываются
\begin{equation}\mathcal{E} _M=\mathcal{E} _H+\mathcal{E} _N.
\end{equation}
Из равенств (2.3) - (2.4) следует, что $$\mathcal{E} _M=\frac{
\dot{\sigma _M}}{2\mu}+\frac{\sigma _M}{2\eta }.
$$

Откуда получаем, что реологическое соотношение для тела Максвелла
имеет вид $$ \sigma _M = e^{ -\frac{\mu }{\eta }t}(\sigma
_{M_0}+2\mu \int\limits_0^te^{ \frac{\mu }{\eta }s}\mathcal{E}
_M(s) \,ds ),$$ где $\sigma_{M_0}$ есть напряжение в начальный
момент времени (здесь и ниже будем считать, что это момент $t=0$).

Одним из типичных представлений моделей вязкоупругих жидкостей
является модель Джеффриса $J=M|N$: параллельно соединяются модель
Максвелла, состоящая из последовательно соединенных пружины и
поршня, и еще один поршень, причем вязкости среды в двух разных
пробирках могут различаться; мы обозначим их $\eta _1$ и $\eta
_2$. Итак, имеем для максвелловской составляющей тела Джеффриса:
\begin{equation}\sigma _M = e^{-\frac{\mu}{\eta _1}t}\left(\sigma
_{M_0}+2\mu \int\limits_0^te^{\frac{\mu}{\eta _1}s} \mathcal{E}
_M(s)\,ds \right),\end{equation} а для ньютоновской составляющей
\begin{equation} \sigma _N = 2\eta _2 \mathcal{E} _N.
\end{equation}
Но при параллельном соединении получаем: \begin{equation}\sigma
_J=\sigma _M+\sigma _N, \end{equation}
\begin{equation}\mathcal{E}
_J=\mathcal{E} _M=\mathcal{E} _N. \end{equation} Равенства
(2.5)-(2.8) влекут
\begin{equation}\sigma _J = e^{-\frac{\mu}{\eta _1}t}\left(\sigma
_{M_0}+2\mu \int\limits_0^t e^{\frac{\mu}{\eta _1}s}\mathcal{E}
_J(s)\,ds \right)+ 2\eta _2\mathcal{E} _J. \end{equation}
Обозначим напряжение и скорость деформации тела Джеффриса в
нулевой момент времени через $\sigma _{J_0}$ и $\mathcal{E}
_{J_0}$. Положив в (2.9) $t=0$, получим равенство:
\begin{equation} \sigma _{J_0}=\sigma _{M_0}+2\eta _2 \mathcal{E}
_{J_0}. \end{equation} Из (2.9) и (2.10) следует $$\sigma _J =
e^{-\frac{\mu}{\eta _1}t}\left(\sigma _{J_0}-2\eta _2 \mathcal{E}
_{J_0}+2\mu \int\limits_0^t \mathcal{E} _J(s)e^{\frac{\mu}{\eta
_1}s}\,ds \right)+2\eta _2\mathcal{E} _J.$$ Это реологическое
соотношение для тела Джеффриса с явно выраженным напряжением.
Более традиционна другая форма реологического соотношения, которую
мы сейчас выведем. Продифференцируем по $t$ выражение (2.9):
\begin{gather}
\dot{\sigma }_J = e^{-\frac{\mu }{\eta _1}t}2\mu e^{ \frac{\mu
}{\eta _1}t} \mathcal{E}_J(t) - \frac{\mu }{\eta _1 }e^{-\frac{\mu
}{\eta _1}t}\left(\sigma_{M_0} + 2\mu \int\limits_0^t
e^{\frac{\mu}{\eta _1}s}\mathcal{E}_J (s) ds\right) + 2\eta_2
\dot{\mathcal{E}}_J.
\end{gather} С другой стороны, из (2.9) следует:
\begin{equation} e^{ - \frac{\mu }{\eta _1}t}\left(\sigma _{M_0} +
2\mu \int\limits _0^te^{  \frac{\mu }{\eta _1}s}\mathcal{E}
_J(s)\,ds \right)=\sigma _J-2\eta _2\mathcal{E} _J.
\end{equation}

Равенства (2.11) и (2.12) влекут
$$ \dot{\sigma }_J=2\mu \mathcal{E} _J - \frac{\mu }{\eta
_1}(\sigma_J -2\eta _2\mathcal{E} _J) + 2\eta _2 \dot{\mathcal{E}
}_J.$$

Умножим это на $ \frac{\eta _1}{\mu }$ и перенесем член $\sigma
_J$ в левую часть:

$$ \sigma _J + \frac{\eta _1}{\mu }\dot{\sigma
}_J=2(\eta _1 + \eta _2)\mathcal{E} _J +2 \frac{\eta _1 \eta
_2}{\mu }\dot{\mathcal{E} }_J.$$

Обозначим $\lambda _1=\frac{\eta _1}{\mu },$ $\lambda _2=
\frac{\eta _1\eta _2}{\mu (\eta _1+\eta _2)},$ $\eta_J=\eta
_1+\eta _2.$ Тогда это перепишется в виде \begin{equation}
\sigma _J + \lambda _1\dot{\sigma }_J=2\eta _J(\mathcal{E} _J
+\lambda _2\dot{\mathcal{E} }_J).\end{equation}

Это и есть другая форма реологического соотношения для тела
Джеффриса, которую мы и будем использовать в дальнейшем.

Уравнение (2.13) является реологическим также и для тела Летерзиха
с символической записью $(N|H)- N$, где опять же вязкости двух
поршней $N$ могут различаться.

\par Далее, применение полученных нами моделей для описания
жидкостей, которые находятся в реальном трехмерном пространстве,
требует осмысления полученных уравнений в ситуации трехмерного
пространства. Отметим, что все обсуждаемые ниже уравнения могут
рассматриваться в пространствах произвольной размерности и тот
факт, что мы до этого момента и ниже говорим о трехмерном
пространстве, объясняется лишь тем, что именно эта ситуация имеет
физический смысл (хотя при исследовании плоскопараллельных течений
или течений в узкой области часто бывает полезно прибегать к
двумерным моделям).

Проблема состоит в том, как понимать составляющие полученных нами
реологических соотношений, например (2.13), в трехмерном случае.
Ответ таков: \emph{времена релаксации} $\lambda_1$ и \emph{запаздывания} $\lambda_2$\footnote{Из вышеприведенных рассуждений следует, что $\lambda_1 > \lambda_2$; это в дальнейшем предполагается по умолчанию.}, а также \emph{вязкость}
тела Джеффриса $\eta$ остаются скалярными величинами, и их можно даже
измерить для конкретных материалов; напряжение и скорость
деформации естественно понимать как тензор касательных напряжений
и тензор скоростей деформации. Остается лишь вопрос о том, как
понимать \emph{производную} по времени, обозначавшуюся точкой.
Однозначного ответа на этот вопрос нет, возможны варианты , и от
выбора такого варианта зависит получаемая модель.

Наиболее простая ситуация, с математической точки зрения,
возникает, если под производной обозначавшейся точкой в (2.13),
понимать частную производную по времени. Получающееся в этом
случае реологическое уравнение \be \sigma+\lambda_1
\frac{\partial\sigma}{\partial t}= 2\eta \left(\mathcal{E}
+\lambda_2 \frac{\partial \mathcal{E}}{\partial t}\right)\ee (здесь
и ниже мы опускаем индекс $J$ в реологическом соотношении)
связывает тензоры касательных напряжений и скоростей деформации и
их изменения в каждой конкретной геометрической точке
пространства.Это существенно ограничивает класс сред, описываемых
данной моделью (см \cite{old2}.)

Более естественным с физической точки зрения представляется
подход, когда реологическое соотношение связывает эти тензоры со
скоростями их изменения для каждой конкретной частицы. Для его
реализации необходимо заменить производную, обозначавшуюся в
(2.13) точкой, на полную или, как еще говорят, субстанциональную
производную $\frac{d}{dt}\equiv \frac {\partial}{\partial t} + v
\cdot \mathrm {grad}$, т.е. в этом случае реологическое
соотношение принимает вид\begin{equation}\sigma+\lambda_1
\frac{d\sigma}{dt}= 2\eta \left(\mathcal{E} +\lambda_2
\frac{d}{dt} \mathcal{E}\right).\end{equation}

Выражая $\sigma$ чисто формально из данного соотношения и
подставляя его выражение в уравнение (0.4), приходим к системе
уравнений следующего вида: $$\frac{\partial v}{\partial t}
+ \sum\limits_{i=1}^n v_i\frac{\partial}{\partial x_i} -
\eta\frac{\lambda_2}{\lambda_1}\triangle v -
2\eta\frac{\lambda_1 - \lambda_2}{\lambda_1^2}\int\limits _0^t
e^{\frac{s-t}{\lambda_1}} \mathrm
{Div}[\mathcal{E}(v)(s,Z(v)(s;t,x))]ds + $$
\begin{equation} + \mathrm {grad} p = f + e^{-\frac{t}{\lambda_1}}
\mathrm {Div}[ \sigma_0 (Z(v)(0;t,x)) - 2 \eta
\frac{\lambda_2}{\lambda_1}\mathcal{E}_0 (Z(v)(0;t,x))],
\end{equation}
$$\mathrm {div} v=0.$$

Уравнения этой системы содержат не только неизвестные скорости
движения $v$ и давления $p$, но и известные траектории
движения частиц среды $z(\tau)=Z(v)(\tau;t,x)$,
определяемые полем скоростей.

Исследование как слабых, так и сильных решений начально-краевой
задачи для уравнений движения (2.16) наталкивается на ту
трудность, что поле скоростей, определяемое решением задачи, не
позволяет восстановить траектории движения частиц или же
траектории не обладают свойством единственности и регулярности,
необходимым для исследования модели.

Однако оказывается, что выражение (2.13) не изменяет форму при
замене переменных (1.1) и (1.2) лишь только тогда, когда под
точкой в реологическом соотношении (2.13) мы будем понимать так
называемую объективную производную. Приведем ее определение. Пусть
$T(t,x)$ -- произвольная тензорнозначная функция, не зависящая от
наблюдателя.

{\bf Определение.}~ Оператор вида
$$\frac{DT(t,x)}{Dt}=\frac{dT(t,x)}{dt}+G(\nabla v(t,x),T(t,x)),
$$ где $G$ -некоторая матричнозначная функция двух матричных
аргументов, называется объективной производной, если при любом
изменении системы отсчета (1.1)-(1.2) выполнено равенство
\begin{equation}\frac{D^*T^*}{Dt^*}=Q(t)\frac{DT}{Dt}Q(t)^T \end{equation}
для всех возможных функций $T$.

\bf Замечание. \rm Символ $\frac{D^*}{D t^*}$ обозначает
представление оператора $\frac{D}{D t}$ в новой системе координат,
т.е. выражение вида
$$\frac{d^* T^*}{dt^*}+G((\nabla v)^*,T^*)=\frac{\partial
T^*}{\partial t^*}+\sum\limits_{i=1}^{3}v^*_i\frac{\partial
T^*}{\partial x^*_i}+G((\nabla v)^*,T^*).$$

Выбор функции $G$ осуществляется на основе различных механических
и опытных соображений.

Имея некоторую объективную производную $\frac{D}{Dt}$, можно от
(2.13)
 перейти к определяющему соотношению
\begin{equation}\sigma +\lambda _1\frac{D\sigma
}{Dt}=2\eta(\mathcal{E}+\lambda
_2\frac{D\mathcal{E}}{Dt}).\end{equation}

Такое определяющее соотношение уже удовлетворяет принципу
объективности. Действительно, подставляя в (2.18) представления
(1.4) и (1.5) и пользуясь свойством (2.17) объективной
производной, получаем:
$$Q(t)^T\sigma ^*Q(t)+\lambda _1Q(t)^T\frac{D^*\sigma
^*}{Dt^*}Q(t)=$$ $$2\eta(Q(t)^T\mathcal{E}^*Q(t)+\lambda
_2Q(t)^T\frac{D^*\mathcal{E}^*}{Dt^*}Q(t)),$$ что влечет $$ \sigma
^*+\lambda _1\frac{D^*\sigma ^*}{Dt^*}=2\eta
\left(\mathcal{E}^*+\lambda
_2\frac{D^*\mathcal{E}^*}{Dt^*}\right).
$$

\par Самый простой пример объективной производной тензора это
производная Яуманна:

$$\frac{D_0T(t,x)}{Dt}=\frac{dT(t,x)}{dt}+T(t,x)W(t,x)-W(t,x)T(t,x).$$

Она также называется  коротационной.

Можно показать \cite{bookour}, что всякая объективная производная может быть
представлена в виде

$$\frac{DT(t,x)}{Dt}=\frac{D_0T(t,x)}{Dt}+G_1(\mathcal{E}(t,x),T(t,x)),$$
где $G_1$ - матричнозначная функция двух матричных аргументов.
Смысл этого представления заключается в том, что любая объективная
производная есть сумма производной Яуманна и выражения, не
зависящего от тензора завихренности $W$.

Простейшее обобщение производной Яуманна - производная Олдройда,
зависящая от параметра $a\in[-1,1]$: \begin{equation}
\frac{D_aT}{Dt}=\frac{D_0T}{Dt} - a(\mathcal{E}T+T\mathcal{E}).
\end{equation}

При $a=0$ это производная Яуманна. При $a=1$ производная (2.19)
называется верхней конвекционной, а при $a=-1$ нижней
конвекционной.

О более сложных моделях типа Джеффриса речь пойдет ниже.

\section{Общий обзор современного состояния математических исследований по моделям типа Джеффриса}

Исторически первым изучался простейший случай, т.е. краевая задача "с частной производной" (0.2),(0.4),(0.5),(2.14), причем последнее соотношение использовалось неявно: $\sigma$ выражалось через $v$ и подставлялось в уравнение движения (0.4), при этом иногда выбрасывались некоторые (по большому счету, несущественные) слагаемые. В этом случае ситуация сходна со случаем Навье-Стокса (т.е. многие результаты, известные для Навье-Стокса, также имеют место), например, в двумерном случае имеется единственное \emph{сильное} решение начально-краевой задачи при всех начальных скоростях, и для них можно построить глобальный аттрактор, а в трехмерном есть (также глобальное по времени) \emph{слабое}  решение, или локальное по времени сильное решение. Из наиболее важных здесь отметим работы \cite{agsu,kos,kos1,osk,nec,ags,ags2,dmz} (последние три посвящены более общему случаю нелинейной вязкоупругости, когда ряд коэффициентов зависит от характеристик скорости, напр. тензора скоростей деформации).

Под сильным и слабым решением в разных работах (мы имеем в виду все работы, описываемые в этом пункте, не только случай с частной производной) понимают несколько разное, но для сильного важно наличие, по крайней мере в соболевском смысле, всех производных, встречающихся в уравнениях. Сильное решение, как правило, единственно, хотя это не всегда строго доказано. У слабого решения скорость имеет лишь первую соболевскую производную по $x$, суммируемую с квадратом. Единственность слабого решения в большинстве рассматриваемых в этом пункте случаев остается открытой проблемой.

Для краевой задачи "с полной производной" (0.2),(0.4),(0.5),(2.15) особенно интересна задача о слабом глобальном решении при заданных начальных скоростях и напряжениях, тогда как специфических результатов о сильных решениях практически нет\footnote{Исключение -- Теорема 5.3 (см. ниже).} (т.е. сильное решение есть тогда, когда оно есть для задач "с объективной" производной, см. ниже, и может быть, очевидным образом, получено тем же методом). Слабая глобальная по времени разрешимость установлена независимо в \cite{turg} и \cite{vz2}. Единственность (при фиксированных плотности внешних сил, начальной скорости и тензоре напряжений) слабого решения не доказана, но имеется \cite{bookour} в классе слабых решений, удовлетворяющих специальному энергетическому неравенству (которому удовлетворяют все сильные решения, а также слабые решения, построенные по методу \cite{vz2}, а возможно (не доказано), что и все слабые решения) при условии, если есть решение из более узкого класса. Cуществует также слабое стационарное (не зависящее от времени) решение данной краевой задачи \cite{v5} (если внешняя сила автономна). Оно принадлежит глобальному аттрактору этой задачи, который был построен в \cite{att2}, используя специально разработанную теорию минимальных траекторных аттракторов (классический подход не проходит из-за отсутствия единственности решения и "инвариантных" диссипативных оценок).  Результаты об аттракторах обобщаются на неавтономный случай \cite{att3,bookour}.

С другой стороны, как мы видели, формально выражая
$\sigma$ из соотношения (2.15), мы приходим к системе
(2.16), в которой необходимо знать траектории движения частиц
среды, определяемые полем скоростей $v(t,x)$. Траектории
движения частиц среды определяются полем скоростей $v$ как
решение интегрального уравнения \begin{equation} z(\tau;t,x) = x+
\int\limits_t^\tau v (s,z(s;t,x))ds.
\end{equation} Если $v$ берется из класса слабых решений, однозначная разрешимость этой системы не имеет места, что ведет к определенным недоразумениям. В связи с этим в \cite{zvd} был предложен вариант
начально-краевой задачи для системы (0.2),(0.4),(0.5),(2.15) со
сглаженным полем скоростей и доказана слабая разрешимость этой
задачи. При стремлении параметра регуляризации к нулю ее решения в
определенном смысле стремятся к решениям нерегуляризованной задачи
(см. \cite{shod}).  Результат о существовании слабого решения
обобщается на движение в областях, меняющихся в зависимости от
времени \cite{orl}. Сильное решение существует локально по времени
\cite{dmizv}; продлить решение глобально можно в случае двумерной
области \cite{dmizv} или при малых внешних силах и начальных
данных \cite{DZ1}.

Перейдем теперь к (нерегуляризованным) моделям с "объективной" производной, т.е. к краевой задаче (0.2),(0.4),(0.5),(2.18). Здесь особняком стоит "коротационный" ("яуманновский") случай \be D=D_0. \ee Дело в том, что в силу кососимметричности матрицы $W$ для всякой симметрической матрицы $A$  той же размерности\footnote{Как обсуждалось выше, основной интерес представляет случай $n=3$.} $n\times n$ имеет место тождество $\sum\limits_{i,j,k=1}^n A_{ij}W_{jk}A_{ik}=0$, что дает возможность сохранить все основные энергетические оценки, имеющиеся для слабых решений задачи (0.2),(0.4),(0.5),(2.15) (ср. с \cite{gran,lmas}). Тем не менее, большую сложность представляет предельный переход в нелинейных членах, содержащих $W$, что не позволяет прямого переноса результатов, имеющихся для задачи (0.2),(0.4),(0.5),(2.15), на коротационную модель. Если же в качестве $D$ взять специальным образом регуляризованную яуманновскую производную (так, чтобы закон (2.18) остался объективным уже с этой регуляризованной производной), то практически все описанные выше (и ниже в п.4) результаты для задачи (0.2),(0.4),(0.5),(2.15)  сохраняются \cite{gran}.

П.Л. Лионсом и Н. Масмуди в статье \cite{lmas} сформулирован результат о глобальной по времени слабой разрешимости при заданных начальных скоростях и напряжениях для задачи (0.2),(0.4),(0.5),(2.18),(3.2). Однако доказательство, представленное ими, имеет крайне сжатый, а местами эвристический характер, и восстановить полное доказательство практически не представляется возможным. В частности, вызывает сомнение следующая аргументация Лионса и Масмуди. На стр. 141 \cite{lmas} для доказательства возможности предельного перехода в последовательности неньютоновских составляющих тензоров напряжений $\tau^n$ вводится функция
$\eta$ такая, что $\|\tau^n\|^2\to \|\tau\|^2 +\eta$
слабо в $L_1(0,T; L_1(\Omega))$, где $\|\cdot\|$ обозначает евклидову норму матрицы. Затем доказывается, что $\eta\equiv 0$, что влечет $\|\tau^n\|^2\to \|\tau\|^2$ слабо в $L_1(0,T; L_1(\Omega))$.
Однако при доказательстве того, что $\eta\equiv 0$, предполагается и не поясняется, что $\eta|_{t=0}= 0$. Но это не является\footnote{Например, если $\tau^n_0\equiv \tau_0\neq 0$ и $\tau^n=\tau_0 sign
\cos nt$, то $\tau^n \to \tau \equiv 0$ слабо в
$L_1(0,T; L_1(\Omega))$,  и $\|\tau^n\|\equiv \|\tau_0\|$, т.е.
$\eta=\|\tau_0\|^2 \neq 0$.} следствием имеющейся сходимости $\tau^n_0\to \tau_0$ в $L_2(\Omega)$.  Остается неясным, верно ли базовое тождество $\eta|_{t=0}= 0$.

Для задачи с производной Яуманна известно (независимо от результата \cite{lmas}), что множество тех правых частей (т.е. внешних сил) начально-краевой задачи, для которых существует сильное глобальное решение, плотно в некоторой топологии \cite{kuz2}, а также получена теорема существования оптимального сильного решения в задаче оптимального управления правыми частями \cite{kuz1}.

В случае более сложных, чем коротационная, производных в задаче
(0.2),(0.4),(0.5),(2.18), использование понятия слабого решения
пока не принесло существенного эффекта, так что
практически все обсуждаемые ниже утверждения касаются сильных решений.  Важно, что
принципиально особых результатов для двумерных областей мало, так
что приведенные ниже верны, как правило, для $n=2,3$.

Для системы с производной Олдройда доказано \cite{gs1} локальное
по времени существование решений начально-краевой задачи для
системы уравнений движения вязкоупругой жидкости в ограниченной
области (результат обобщен в \cite{tal} на неограниченные области,
а в \cite{tal1} на случай более общих, чем (0.5), граничных
условий). Глобальное существование имеется при малых начальных
данных и внешних силах. Впервые это показано в \cite{gs1} в
ограниченной области при дополнительном условии, что времена
релаксации $\lambda_1$ и запаздывания $\lambda_2$ близки; в случае
$\Omega=\R^n$ это условие лишнее \cite{bookour,cm}, а затем
оказалось \cite{mol1}, что от него можно избавиться и для
ограниченной области; подобный результат, полученный другим
методом и адаптированный для численных схем, имеется в
\cite{bonito} при пренебрежении квадратичными инерционными
(конвективными) членами в (2.18) и (0.4). В \cite{lei} схожие
локальные и глобальные (при малых данных) теоремы получены для
сжимаемой жидкости, а затем, путем предельного перехода по
параметру, отвечающему за сжимаемость, и для несжимаемой.
Указанные утверждения (кроме \cite{bonito}) относятся к
рассмотрению решений в гильбертовых соболевских пространствах типа
$W^r_2$ по $x$. Аналогичные (иногда с дополнительными
ограничениями) результаты получены для банаховых пространств типа
$W^r_p$ \cite{fer,fer1}. В \cite{cm} изучался случай $\Omega=\R^n$
(или тор, в обоих случаях граничное условие (0.5) опускается) и
также получены теоремы подобного содержания в пространствах
Бесова, особое внимание уделено классам единственности. При малой
периодической по времени правой части задача имеет периодическое  решение \cite{mol1}. При малой же
не зависящей от времени правой части задача имеет стационарное
решение \cite{ren,tal3} (затем результат обобщался на различные
неограниченные области \cite{pil,gt}). Устойчивость стационарного решения
изучалась в \cite{renar1,gs3} (см. также \cite{renar} с обзором по
поводу устойчивости и неустойчивости систем уравнений этого и
более широких классов). При стремлении к нулю времени релаксации
(а значит, и не превосходящего его времени запаздывания) решения
(поля скоростей) стремятся к решениям задачи Навье-Стокса
\cite{mol2} (здесь предполагается, что $\Omega=\R^n$ или тор),
причем при малых временах релаксации все решения существуют на том
же интервале, что и решение для Навье-Стокса. Отсюда выводится
важное  следствие \cite{mol2} о глобальном по времени
существовании сильного решения для $\Omega=\R^2$ (или квадрата с
периодическими граничными условиями), произвольных данных и малом
времени релаксации (хотя степень его ''малости'', разумеется,
зависит от остальных данных и параметров). В \cite{lin} и других
работах этой группы ученых был предложен другой подход к данной
проблеме, в т.ч. механически мотивированные замены переменных, а
также добавки и изменения членов, дающие возможность улучшить
математические свойства задачи. По поводу численных схем для
начально-краевой задачи для системы уравнений движения
вязкоупругой жидкости с производной Олдройда см. напр.
\cite{bonito,boy}.

Математическое исследование более общих задач, чем (0.2),(0.4),(0.5),(2.18) с производной Олдройда, можно условно разбить на три группы.

Первая группа связана с обобщением определяющего соотношения (2.18) с производной Олдройда. Самое простое обобщение --- предположение, что вязкость зависит от евклидовой нормы тензора скоростей деформации (т.е. его второго инварианта). Существование стационарного решения краевой задачи при малой автономной внешней силе в разных типах областей показано в \cite{arad,arad1}. Другое обобщение --- модель Уайта-Метцнера, где уже время релаксации $\lambda_1$ и параметр вязкости
$(1-\lambda_2/\lambda_1)\eta$ зависят от второго инварианта тензора скоростей деформации.
Краевые задачи для этой модели, во многом, имеют свойства, похожие на задачи с производной Олдройда. Основные результаты, получаемые здесь, изложены в \cite{ght,hak,mol0}. Еще одно обобщение --- комбинированная модель для смеси нескольких вязкоупругих (с объективной производной в общем виде) и нелинейно-вязких жидкостей \cite{bookour}. Соответствующая начальная задача c
$\Omega=\R^n$ имеет глобальное решение при малых данных \cite{vz3,bookour}.

Вторая группа содержит задачи c реологическим законом типа (2.18) с производной Олдройда, где имеется дополнительная переменная, от которой зависят вязкость, времена релаксации и запаздывания. В качестве такой переменной может выступать температура (в этом случае в \cite{besbes} показано локальное существование решений) или процентное соотношение в смеси двух вязкоупругих жидкостей (в \cite{chup} имеются утверждения о локальном существовании решений, а также глобальном при малых данных).

Третья группа --- это так называемые микро-макро модели, учитывающие соображения из квантовой механики. Математически, однако, соответствующие начально-краевые задачи (см. напр. \cite{barr,otto,lmas1}) обладают лучшими свойствами с точки зрения диссипативности, что облегчает вопросы разрешимости и предельного поведения, несмотря на громоздкость систем уравнений.

\section{О слабых решениях начально-краевых задач для модели Джеффриса
с субстанциональной производной}

В этом параграфе будут описаны исследования по проблеме
разрешимости начально-краевой задачи (0.2),(0.4),(0.5),(2.15)\footnote{плюс начальное условие}в слабом
смысле.

Вначале введем используемые функциональные пространства.

Пусть $E$ --- какое-то конечномерное вещественное линейное
пространство, например, $\mathbb{R}^n$. Мы будем использовать
стандартные обозначения $L_p(\Omega,E),~
W_p^m(\Omega,E),~H^m(\Omega,E)=W_2^m(\Omega,E),
~H_0^m(\Omega,E)=\stackrel{\circ}{W^m_2}(\Omega,E)$ пространств
Лебега и Соболева для функции из области
$\Omega\subseteq\mathbb{R}^n $ со значениями в $E$. Иногда для
краткости будем писать просто $L_p$ вместо $L_p(\Omega,E)$ и т.п.,
если из контекста ясно, о каком пространстве $E$ идет речь.

Скалярные произведения в $L_2(\Omega,E)$  и $H^1(\Omega,E)$ могут
быть заданы выражениями $$ (u,v) = \int\limits _\Omega
(u(x),v(x))_Edx,$$ $$(u,v)_1 = (u,v)+
\sum\limits_{i=1}^{n}\left( \frac{\partial u}{\partial
x_i},\frac{\partial v}{\partial x_i} \right ).$$

Евклидову норму в $E$ будем обозначать символом $|\cdot|$, норму в
$L_2$ - символом $\| \cdot \|$, а норму в $W_2^1$ - символом
$\|\cdot\|_1$.

Через $C_0^\infty(\Omega,E)$ будем обозначать пространство гладких
функций с компактным носителем в $\Omega$ и со значениями в $E$.

Пусть $\mathcal{V}=\{u\in C_0^\infty(\Omega,\mathbb{R}^n),\mathrm
{div} u = 0 \}$.

Символами $H$ и $V$ обозначаются замыкания $\mathcal{V}$
соответственно в $L_2(\Omega,\mathbb{R}^n)$ и
$W_2^1(\Omega,\mathbb{R}^n)$.

Через $H^*$ и $V^*$ будем обозначать сопряженные пространства к
соответственно пространствам $H$ и $V$.

Как обычно делается в математической гидродинамике \cite{bookour}, будем отождествлять пространство $H$ и его
сопряженное пространство $H^*$. Поэтому имеем вложение $$V\subset
H \equiv H^*\subset V^*.$$

Символом $\langle f,v \rangle$ будем обозначать действие
функционала $f\in V^*$ на элемент $v\in V.$ В силу указанного выше
отождествления, если $f,v\in V$, то $\langle f,v \rangle = (f,v).$

Символами $C([0,T];X), C_\omega([0,T];X), L_2(0,T;X)$ и т.п.
обозначаются банаховы пространства непрерывных, слабо непрерывных,
суммируемых с квадратом и т.п. функций на промежутке $[0,T]$ со
значениями в некотором банаховом пространстве $X$.

Далее, как  уже было отмечено в предыдущем параграфе, выражая
$\sigma$ чисто формально из соотношения (2.15), приходим к системе
(2.16), в которой необходимо знать траектории движения частиц
среды, определяемые полем скоростей $v(t,x)$. Траектории движения
частиц среды определяются полем скоростей $v$ как решение
интегрального уравнения \begin{equation} z(\tau;t,x) = x+
\int\limits_t^\tau v (s,z(s;t,x))ds.
\end{equation}
Но, как мы увидим ниже, слабые решения $v (t,x)$
начально-краевой задачи для системы (0.2), (0.4),(2.15)принадлежат
пространству \begin{equation} W_{2,1}= \{ v;v \in
L_2(0,T;V),v' \in L_1(0,T;V^* ) \}
\end{equation} и не являются непрерывными функциями по переменной
$x$ для $n=2,3.$

Однако существование решений уравнения (4.1) при фиксированном
$v$ известно лишь для случая $v \in
L_1(0,T;C(\Omega)^n)$ (cм. \cite{OS2}) и это решение единственно для
$v \in L_1(0,T;C^1(\overline{\Omega})^n)$ (см. там же). В
связи с этим фактом в \cite{zvd} был рассмотрен регуляризованный вариант
начально-краевой задачи для системы (0.2),(0.4),(2.15). Опишем
его.

Выберем некоторый линейный оператор регуляризации $S_\delta : H
\rightarrow C^1(\overline{\Omega})^n \cap V$ для $\delta > 0$
такой, что порождаемое им отображение $S_\delta : L_2(0,T;H)
\rightarrow L_2(0,T;C^1(\overline{\Omega})^n \cap V)$ непрерывно,
а операторы $S_\delta : L_2(0,T;H) \rightarrow L_2(0,T;H)$
сходятся сильно к тождественному оператору $I$ при $\delta
\rightarrow 0$. Конструкция такого оператора приведена в
\cite{zvd}, уточнена в \cite{DZ1}и подробно описана в
\cite{regul}.

Обозначим через $C^1 D (\overline{\Omega})$ множество взаимно
однозначных отображений $z:\overline{\Omega} \rightarrow
\overline{\Omega}$, совпадающих с тождественным отображением на
$\partial\Omega$  и имеющих непрерывные частные производные
первого порядка такие, что $det \left( \frac{\partial z}{\partial
x} \right)= 1$ в каждой точке области $\overline{\Omega}$. Будем
предполагать, что в этом множестве используется метрика
пространства непрерывных функций
$C(\overline{\Omega},\mathbb{R}^n)$.

Положим также $CG = C([0,T]\times [0,T], C^1 D
(\overline{\Omega}))$.

Рассмотрим уравнение $$ z(\tau ;t,x)= x + \int \limits_t^\tau
S_\delta v(s,z(s;t,x))ds ,~~~~~ \tau ,t \in (0,T), { }{
}x\in\Omega.
$$

Для каждого $v\in L_2(0,T;V)$ это уравнение имеет
единственное решение в классе $CG$ (см. \cite{OS2}).Обозначим его через
$Z_\delta (v)(\tau;t,x).$

Полагая для простоты в (2.16) $$\sigma_0 (Z_\delta
(v)(0;t,x)) - 2\eta
\frac{\lambda_2}{\lambda_1}\mathcal{E}_0
(Z_\delta(v)(0;t,x)) = 0,
$$ мы приходим к следующей начально-краевой задаче для
регуляризованной модели, описывающей движение несжимаемой
вязкоупругой среды: $$\partial_tv (t,x) +
\sum\limits_{j=1}^n v_j \frac {\partial v}{\partial
x_j}(t,x)-\mu_1 \mathrm {Div}\int\limits_0^t
e^{-\frac{t-s}{\lambda_1}} \mathcal{E}(v)(s,z(s;t,x))ds -$$
\begin{equation} - \mu_0 \Delta v(t,x)+
 \mathrm {grad}p(t,x) =
f(t,x), ~~~~~~~{}(t,x)\in(0,T)\times\Omega,
\end{equation}
\begin{equation} z(s;t,x)=x +\int\limits_t^s S_\delta
v (\tau,z(\tau;t,x))d\tau, ~~~~{} s\in[0,T],~~~~
(t,x)\in(0,T)\times\Omega,
\end{equation}
\begin{equation}\mathrm {div}v = 0,~~~~~~~~
(t,x)\in(0,T)\times\Omega,
\end{equation}
\begin{equation} v |_{(0,T)\times\partial\Omega}=0,
\end{equation}
\begin{equation}v(0,x)=v^0(x),~~~~~~{}x\in\Omega, \end{equation}
\begin{equation} \int\limits_\Omega p{}~dx = 0.
\end{equation}

Здесь $\mu_0 = \eta \frac{\lambda_2}{\lambda_1}$,~~$\mu_1 = 2\eta
\frac {\lambda_1 - \lambda_2}{\lambda_1^2}.$

Как уже отмечалось, при введении оператора регуляризации уравнение
(4.4) имеет единственное решение
$Z_\delta(v)(s;t,x).$Подставив его в (4.3), получим
следующую начально-краевую задачу, эквивалентную задаче
(4.3)-(4.8):
$$\partial_t v (t,x) + \sum
\limits_{j=1}^n v_j \frac {\partial v}{\partial
x_j}(t,x) -$$
\begin{equation} -\mu_0 \Delta v(t,x) -\mu_1
\mathrm {Div} \int\limits_0^t
e^{-\frac{t-s}{\lambda}}\mathcal{E}(v(s,Z_\delta(v)
(s;t,x))ds + \end{equation}
$$+ \mathrm {grad}p(t,x)=f(t,x), ~ ~~~~~~~~~~ (t,x)\in (0,T)\times
\Omega $$
\begin{equation}\mathrm {div} v(t,x)= 0,~~~~~~~~(t,x)\in (0,T)\times
\Omega ,\end{equation}
\begin{equation} v(t,x) = 0, ~~~~~~~(t,x)\in(0,T)\times
\partial \Omega, \end{equation}
\begin{equation}v(0,x) = v^0(x),~~~~~~~~x\in\Omega, \end{equation}
\begin{equation}\int \limits_\Omega^{} p(t,x) dx = 0. \end{equation}

 {\bf Определение 4.1.}~~Пусть $f\in L_1(0,T;V^*).$ Слабым решением начально - краевой задачи (4.9)-(4.13)
 называется вектор-функция $v\in W_{2,1}$,
 удовлетворяющая начальному условию (4.12)и интегральному тождеству
 $$ \frac{d}{dt}(v(t),\varphi)-
 \sum\limits_{i=1}^{n} \left( v_i (t) v (t) , \frac{\partial\varphi}{\partial x_i}\right) + $$
 \begin{equation}+ \mu_0(\nabla v(t),\nabla \varphi)+\mu_1\left(\int\limits_0^t
 e^{\frac{s-t}{\lambda_1}}\mathcal{E}(s,Z_\delta(v)(s;t,x)ds,\mathcal{E}(\varphi)\right)=\langle f(t),\varphi\rangle \end{equation}
 для всех $\varphi \in V$ п.в. на $t\in(0,T)$.

Используя аппроксимационно-топологический метод, изложенный в
\cite{zvdkn,bookour}, в работе \cite{zvd} доказана следующая

{\bf Теорема 4.1.} ~Для каждой функции $f\in L_1(0,T;H^*)+
L_2(0,T;V^*)$ и $v^0 \in H$ задача (4.9)-(4.13)имеет хотя
бы одно слабое решение $v \in W_{2,1}$.

В работе \cite{vz2} рассмотрена начально - краевая задача для исходной
модели Джеффриса без регуляризации поля скоростей. В этом случае
нельзя явно выразить тензор $\sigma$ через вектор-функцию
$v$ и приходится рассматривать систему с переменными
$v = (v_i), p$ и $\sigma = (\sigma_{ij})$.

Итак, запишем начально - краевую задачу для системы
(0.2),(0.4),(2.15) в виде \begin{equation} \frac{\partial
v}{\partial t} +\sum\limits_{i=1}^{n} v_i
\frac{\partial v}{\partial x_i} +\mathrm {grad} p - \mathrm
{Div} \sigma = f(t,x), ~~~t\in (0,T),~x\in\Omega,  \end{equation}
\begin{equation} \sigma + \lambda_1 \left( \frac{\partial \sigma}{\partial
t} + \sum\limits_{i=1}^{n} v_i \frac{\partial
\sigma}{\partial x_i}\right ) $$ $$= 2\eta \left(\mathcal{E} +
\lambda_2 \left(\frac{\partial \mathcal{E} }{\partial t} +
\sum\limits_{i=1}^{n} v_i \frac{\partial
\mathcal{E}}{\partial x_i}\right) \right),~~~t\in (0,T),~
x\in\Omega,\end{equation}
\begin{equation} \mathrm
{div} v (t,x)= 0,~~~t\in (0,T),~x\in\Omega \end{equation}
\begin{equation} v |_{\partial \Omega} = 0,~~t\in (0,T),
\end{equation}
\begin{equation} v |_{t= 0}=v^0,~~\sigma |_{t= 0}=
\sigma_0. \end{equation}

Обозначим через $\R^{n\times n}$ пространство вещественных матриц
порядка $n\times n$, а через $\R^{n\times n}_S$ - его
подпространство симметрических матриц.

{\bf Определение 4.2.}~Пусть $f\in L_1(0,T;V^*).$ Слабым решением
задачи (4.15)-(4.19) называется пара $(u,\sigma)$,
\begin{equation}u\in L_2(0,T;V) \cap C_w([0,T];H),~~~~\frac{du}{dt}\in L_1(0,T;V^*), \end{equation}
\begin{equation}\sigma \in L_2(0,T;L_2(\Omega, \R^{n\times n}_S)) \cap C_w([0,T];H^{-1}(\Omega,
\R^{n\times n}_S)),\end{equation} удовлетворяющая условию (4.19)и тождествам
$$ \frac{d}{dt}(u,\varphi) + (\sigma,\nabla \varphi)-
\sum\limits_{i=1}^{n}\left(u_i u, \frac{\partial \varphi}{\partial
x_i}\right) = \langle f,\varphi \rangle,$$
$$(\sigma , \Phi) + \lambda_1 \frac {d}{dt}(\sigma , \Phi)-\lambda_1
\sum\limits_{i=1}^{n}\left(u_i \sigma, \frac{\partial
\Phi}{\partial x_i}\right)=$$
$$= 2\eta (\mathcal{E}(u),\Phi) + 2\eta \lambda_2 \left( \frac {d}{dt}
(\mathcal{E}(u),\Phi)+ \sum\limits_{i=1}^{n}\left(u_i \mathcal
{E}(u), \frac{\partial \Phi}{\partial x_i}\right)\right)$$в смысле
распределений на $(0,T)$ для всех $\varphi \in \mathcal{V}$ и
$\Phi\in C_0^\infty(\Omega,\R^{n\times n}_S).$

В \cite{vz2} была доказана следующая

{\bf Теорема 4.2.} Пусть $f\in L_2(0,T;V^*),~ v^0\in H,
\sigma_0 \in W_2^{-1}(\Omega,\R^{n\times n}_S),~ \sigma_0-2\eta
\frac{\lambda_2}{\lambda_1}\mathcal {E}(v^0)\in
L_2(\Omega,\R^{n\times n}_S).$Тогда существует слабое решение задачи
(4.15)-(4.19)в классе (4.20)-(4.21).
\par
Конечно, естественно возникает вопрос, как связаны между собой
решения начально-краевой задачи для исходной модели Джеффриса и
решения для регуляризованного варианта (4.3)-(4.8)~этой модели.

Для ответа на этот вопрос будем предполагать, что
$v^0~\in~H,~\sigma_0\in W_2^{-1}(\Omega,\R^{n\times n}_S),~ f\in
L_2(0,T;V^*) $ и выполнено условие
$$\sigma_0(t,x) = 2\eta \frac{\lambda_2}{\lambda_1}\mathcal
{E}(v^0)(t,x).$$

Из теоремы 4.1 следует, что для каждого $\delta > 0$ задача
(4.3)-(4.8) при этих условиях имеет решение $v_\delta$. Это
решение порождает соответствующий тензор касательных напряжений
\begin{equation}\sigma(u)(t,x)= 2\mu_0\mathcal {E}(u)(t,x) + \mu_1
\int\limits_{0}^{t} e^{-\frac {t-s}{\lambda}}
\mathcal{E}(u)(s,Z_\delta (u)(s;t,x))ds. \end{equation}

Имеет место следующая

{\bf Теорема 4.3} (см. \cite{shod}).~Пусть $ \{ \delta_k \}$ -
некоторая последовательность положительных чисел, стремящаяся к
нулю. Пусть $v_k$ - решение задачи (4.9)-(4.13), соответствующее
$\delta_k$, и $\sigma_k$ - соответствующий тензор касательных
напряжений, найденный по формуле (4.22).

Тогда существует пара $(v_*,\sigma_*)$, являющаяся слабым
решением задачи (4.15)-(4.19), для которой существуют
подпоследовательности $\{ v_{k_m} \}$ и $\{\sigma_{k_m}\}$
такие, что

1) $\{v_{k_m}\}$ сходится слабо к $v_*$ в
$L_2(0,T;V);$

2) $\{v_{k_m}\}$ сходится сильно к $v_*$ в
$L_2(0,T;H);$

3)  $\{v_{k_m}\}$ сходится $*$ - слабо к $v_*$ в
$L_\infty(0,T;H);$

4)  $\{\sigma_{k_m}\}$ сходится слабо к $\sigma_*$ в
$L_2(0,T;L_2);$

5) $\{\sigma_{k_m}\}$ сходится $*$ - слабо к $\sigma_*$ в
$L_\infty(0,T;H^{-1}).$

В работе \cite{orl} исследована разрешимость начально-краевой задачи
для регуляризованной модели Джеффриса в случае движения
вязкоупругой среды в переменной области (области, меняющейся в
зависимости от времени).

Опишем результат этой работы более подробно.

Пусть $\Omega_{t}\in R^n,\ n\ge 2$, является семейством
ограниченных областей с границей $\Gamma_t$, $Q=\{ (t,x): t\in [0,
T],\ x\in \Omega_t\}$, $\Gamma =   \{ (t,x): t\in[0, T],\ x\in
\Gamma_t\}$. Рассматривается    начально-краевая задача $$
v_{t} + \sum \limits_{i=1}^{n} \frac {v_i
\partial v}{\partial {x_i}}-2\mu_0  \mathrm {Div} \mathcal{E}(v)
- \mu_1 \mathrm {Div} \int_0^t e^{-\frac {t - s}{\lambda}}
\mathcal{E}(v)(s,z(s;t,x)ds + $$
\begin{equation}
 + \mathrm{grad} p = \varphi,~~\mathrm{div}v =0,~ (t,x) \in Q;
~\int_{\Omega _t } p\,dx = 0,\ t \in [0, T];
\end{equation}
 $$ v(0,x) = v^0 (x), \  x\in \Omega_0,\   v(t,x)=v^1 (t,x)\in
\Gamma.$$ 

Здесь, как и выше, вектор-функция $v(t,x) = (v_1,..., v_n)(t,x)$
--- скорость среды в точке $x$ в момент $t$, скалярная функция
$p$ --- давление, $\mu_0$, $\mu_1$, $\lambda$- положительные
константы. Функция $z(s;t,x)$ определяется как решение задачи Коши
\begin{equation}
z(s;t,x) = x + \int_t^{s}v(\tau,z(\tau;t,x))\, d\tau,\ \ \ s \in
[0,T],\ \ (t,x) \in Q.
\end{equation}
 Предполагается, что $Q\subset R^{n+1}$
 определяется  эволюцией семейства $\Omega_{t}$, $t\ge 0$ объема
 $\Omega_{0}$  с достаточно гладкой границей $\Gamma_0$ вдоль поля
 скоростей некоторого  достаточно гладкого соленоидального векторного
 поля $\tilde v(t,x)$,
 определенного в некоторой цилиндрической области $\hat Q_0=\{ (t,x):
t\in[0, T]$,$\ x\in \hat \Omega_0\}$, так что $\Omega _t$$\subset
\hat \Omega_0$. Это означает,  что $\Omega_{t}=$$\tilde
z(t;0,\Omega_0)$, где $\tilde z(\tau;t,x)$  решение задачи Коши
\begin{equation}
\tilde z(\tau;t,x) = x + \int_t^{\tau} \tilde v(s,\tilde
z(s;t,x))\, ds,\ \ \ \tau \in [0,T],\ \ (t,x) \in \hat{Q_0}.
\end{equation}

Введем необходимые функциональные пространства.
 Обозначим нормы в $L_2 (\Omega_t)$ и $W_2^k(\Omega_t)$  через
$|\cdot |_{0,t}$  и  $| \cdot |_{k,t}$ соответственно, а через
$||\cdot ||_{0}$ обозначим норму в $L_2 (Q)$. Обозначим через
$D_{0,t}$ множество гладких соленоидальных
 финитных  функций в области $\Omega_t$. Обозначим
через $H_{t}$ и $V_{t}$ замыкание $D_{0,t}$ в нормах  $L_2
(\Omega_t)$ и $W_2^1(\Omega_t)$ соответственно.  Сопряженное к
$V_{t}$ пространство обозначим через $V_{t}^{*}$, норму в
$V_{t}^{*}$ обозначим через $|\cdot |_{-1,t}$, а действие
функционала  $v\in V_{t}^{*}$ на элементе $h\in V_t$ обозначим
через $\langle v,h\rangle_t$. При этом скалярное произведение
$(\cdot,\cdot)_t$ в $H_{t}$ порождает
 при каждом $t\in [0,T]$ плотное непрерывное вложение
$ V_{t}\subset H_{t}\subset V_{t}^{*}. $

 Обозначим через  $D$ множество гладких вектор-функций на $Q$,
 соленоидальных и финитных   в области $\Omega_t$ при каждом $t$.
  Обозначим через
  $E$, $E^*$, $E_1^*$,  $W$, $W_1$, $CH$, $EC$, $LH$, $L_{2,\sigma}(Q)$
 замыкание $D$ соответственно в нормах
$$
 \| v \|_{E}\! =\!( \int_0^{T}\!| v(t,x)|^2_{1,t}\,dt)^{1/2},$$ $$
 \|v\|_{E^*}\! =\!(\int_0^{T}\!| v(t,x)|^2_{-1,t}\,dt)^{1/2},\
 \|v\|_{E_1^*}=\!
 \int_0^{T}\!|v(t,x)|_{-1,t}\,dt,
$$
$$
 \|v\|_{W} =\|v\|_{E}+\|v_t\|_{E^*},\
 \|v\|_{W_1} =\|v\|_{E}+\|v_t\|_{E_1^*},\
\|v\|_{CH}= \max\limits _{t\in [0,T]}|v(t,x)|_{0,t},\
$$
$$
\|v\|_{EC}=\|v\|_{E}+\|v\|_{CH},\
\|v\|_{LH}=\int_{0}^{T}|v(t,x)|_{0,t}\,dt,$$ $$
 \|v\|_{0}=(\int_0^{T}|v(t,x)|^2_{0,t}\,dt)^{1/2}.
$$

Пространство  $E^*$ является сопряженным к пространству $E$, а
скалярное произведение $(v,u)=\int_{0}^{T}(v(t,x),u(t,x))_{t}\,dt$
в $L_{2,\sigma}(Q)$ порождает плотное непрерывное вложение $
\label{f21} E  \subset L_{2,\sigma}(Q) \subset E^*$. Пространство
$E^*_1$ непрерывно вложено в $E^*$.

Ниже через $\langle v,u\rangle$ обозначается  действие функционала
$v$ из $E^*$ на функцию $u$ из $E$.

Задача (4.23) включает интеграл, вычисляемый вдоль траектории
$z(s;t,x)$ движения частиц в поле скоростей $v(t,x)$.
Однако, как и в случае цилиндрических областей, даже сильные
решения $v(t,x)\in W_2^{1,2}(Q)$ задачи (4.23) для $n=2,3$
не обеспечивают однозначную разрешимость задачи (4.24). В качестве
выхода из этой ситуации, как и выше, мы регуляризуем
поле скоростей с помощью введения сглаживающего 
оператора $ S_{\delta ,t}: H_t \to C^1( \overline{\Omega}_t) \cap
V_t$ для $\delta >0$ такого, что $S_{\delta ,t} (v)$ близко к $v$.
Определим на $L_{2,\sigma}(Q)$ сглаживающий  оператор
$\hat{S}_\delta (v)=\hat{v} $, где $\hat{v}(t,x)=S_{\delta
,t}(v(t,x))$ при $t \geq 0.$ Оператор
$\hat{S}_\delta:L_{2,\sigma}(Q)\rightarrow L_2(0,T;C^1)\cap E$
оказывается ограниченным.
Пусть теперь $u(t,x)\in L_{2,\sigma}(Q)$,
$v(t,x)=u(t,x)+\tilde{v}(t,x)$. 
      Определим оператор регуляризации $S_\delta : L_{2,\sigma }(Q) \to
L_2(0,T;C^1)\cap E$ формулой $ S_\delta (v)= \hat{S}_\delta
(v-\tilde{v}) +\tilde v = \hat{S}_\delta (u) +\tilde v$.

Заменим (4.24)на уравнение
\begin{equation}
z(\tau; t,x) = x+ \int_{t}^{\tau} S_\delta v(s,z(s;t,x))\, ds, \;
\tau ,t \in [0,T], x \in \Omega _t. \end{equation} Для каждой
$v(t,x)\in E$ функция $S_\delta (v) \in L_2(0,T;C^1),$ и,
следовательно, задача (4.26) однозначно разрешима. Решение задачи
(4.26) будем обозначать $\tilde Z_{\delta} (v)$.


Нас будет интересовать вопрос о разрешимости следующей
регуляризованной задачи:
$$ v_t +  \sum \limits_{k=1}^{n}v_k
\frac {dv}{dx_k} - 2\mu_0 \mathrm{Div}
\mathcal{E}(v) - \mu_1 \mathrm{Div} \int_{0}^{t}
e^{-\frac{t-s}{\lambda}} \mathcal{E}(v)
(s,Z_{\delta}(v)(s;t,x))ds + \nabla p = \Phi ,$$
\begin{equation}
  \mathrm{div} \,v(t,x)=0,\ (t,x)\in Q;  \int_{\Omega _t}{} p(t,x)\,dx=0,\
t\in [0,T]; \end{equation}
$$
v(0,x)= v^0(x), \; x\in \Omega _0, \; v(t,x)=v^1(t,x), \; (t,x)\in
\Gamma, \ t\in [0,T].$$

{\bf Определение 4.3.}~~Пусть $ \Phi (t,x)\in V_t^* $ при п.в.
$t$. Слабым решением задачи (4.27) называется функция $v(t,x)$
вида $v(t,x)=\tilde v (t,x)+w(t,x)$, $w(t,x)\in
E\bigcap CH$, $w(0,x)=0$, такая что при 
 любых $h(t,x) \in D$, таких, что $h(T,x)=0$, справедливо тождество
$$-\int_0^T(v(t,x),h_t(t,x))_t\,dt +\mu_0\int_0^T(v_i(t,x)v_j(t,x),\partial
h_i(t,x)/\partial x_j )_t\,dt +
$$
\begin{equation}
 +\mu_1\int_0^T( \int_0^{t}
e^{-\frac{(t-s)}{\lambda}}\mathcal{E}(v)(s,Z_{\delta
}(v)(s;t,x))\,ds, \mathcal{E}(h(t,x)))_t\,dt =\end{equation}
$$ \langle \Phi (t,x),h(t,x)\rangle - \langle \tilde{v}(0,x),h(0,x)\rangle_0.
$$
Сформулируем  основной результат.

{\bf Теорема 4.4} \cite{orl}. Пусть $2\le n\le 4$, $\Phi =f_1 +f_2,
f_1\in E^*_1$, $f_2\in E^*,$ $v^0(x)=\tilde v(0,x), x\in \Omega_0,
v^1 (t,x)= \tilde v(t,x), (t,x)\in \Gamma.$
Тогда задача (4.27) 
  имеет по крайней мере одно слабое решение.

{\bf Замечание.} В заключении этого параграфа отметим, что после
доказательства разрешимости в слабом смысле задачи (4.15)-(4.19),
может возникнуть вопрос о необходимости исследования разрешимости
в слабом смысле регуляризованого варианта (4.9)-(4.13) этой
задачи. Один из аргументов в пользу рассмотрения регуляризованной
задачи состоит в том, что она представляет интерес с точки зрения
нахождения численных решений для модели Джеффриса, и для этого,
конечно, необходимо знать в каких функциональных пространствах
имеет место разрешимость соответствующей начально - краевой
задачи.

\section{О сильных решениях начально-краевой задачи для модели Джеффриса
с субстанциональной производной}

В этом параграфе будут изложены результаты о существовании сильных
решений начально - краевой задачи для регуляризованной модели
Джеффриса.

Пусть $\Omega$ - ограниченная область в $\mathbb{R}^n$ с границей
$\partial \Omega$ класса $\mathcal{C}^2$. Через $Q_T = [0,T]\times
\Omega$ обозначается пространственно-временной цилиндр для $T>0$ и
через $(t,x)$ -  точки $Q_T$.

Рассмотрим движение среды, заполняющей область $\Omega$ в
$\mathbb{R}^n$, $n = 2,3,$ на промежутке времени $(0,T),~~T>0.$
Ниже будем предполагать, что движение этой среды описывается
системой уравнений
$$\rho \left(\frac {\partial v}{\partial t} + \sum
\limits_{i=1}^{n} v_i \frac {\partial v}{\partial
x_i}\right)(t,x)- \mu_1 \mathrm {Div}\int \limits_0^t
e^{-\frac{t-s}{\lambda}}\mathcal{E}(v)(s,z(s;t,x))ds -
2\mu_0 \mathrm {Div}\mathcal{E}(v)(t,x)+$$
\begin{equation}+ \mathrm {grad}p (t,x) = e^{-\frac{t}{\lambda}}\mathrm
{Div}\sigma_0(z(0;t,x))+
\rho\varphi(t,x),~~~~~~~~(t,x)\in(0,T)\times \Omega,\end{equation}
\begin{equation}z(\tau;t,x)=x + \int\limits_{t}^{\tau}S_\delta
v(s,z(s;t,x))ds,~~~~~~~\tau\in[0,T]\times\Omega,
\end{equation}
\begin{equation}\mathrm {div}v = 0,
~~~~(t,x)\in(0,T)\times\Omega,
\end{equation} где, как и в предыдущих параграфах, $v =
(v_1,\ldots,v_n)$ - скорость движения среды, $p$ -
давление, $\varphi$ - плотность внешних сил, $\rho = const$ -
плотность среды, $\sigma_0$ - тензор остаточных напряжений при
$t=0$, $\mathcal{E} = (\varepsilon_{ij})$ - тензор скоростей
деформации и  $S_\delta$ оператор, определенный в $\S~4.$ Будем
предполагать также, что заданы начально- краевые условия
\begin{equation}v(t,x)=0~~~~~~~~(t,x)\in(0,T)\times
\partial\Omega, \end{equation}
\begin{equation}v(x) = v^0 (x),~~~~~~~x\in\Omega.
\end{equation}

Из уравнений $(5.1)$ видно, что вместе с функцией $p$ этому
уравнению удовлетворяют и все функции, отличающиеся от $p$ на
константу. Поэтому для определенности в постановке
начально-краевой задачи для системы $(5.1)-(5.3)$добавляют условие
\begin{equation} \int\limits_{\Omega}^{}p~dx = 0.
\end{equation}

Введем обозначение $${H_s}^{1,2}(Q_T)^n =
L_2(0,T;W_2^2(\Omega)^n\cap V)\cap L_2(\Omega ; W_2^1(0,T)^n).$$

Теперь мы готовы дать следующее

{\bf Определение 5.1.} Пусть $\sigma_0\in W_2^1(\Omega)^{n^2},
\varphi \in L_2(0,T;L_2(\Omega)^n), v^0 \in V$ - заданные функции
и $n=2,3.$ Сильным решением задачи $(5.1)-(5.6)$ называется пара
$(v, p)$, где $v \in H_s^{1,2}(Q_T)^n ,~p \in L_2(0,T;
W_2^1(\Omega)),$ удовлетворяющие уравнениям (5.1), (5.3) почти
всюду в $Q_T$, начальному условию $(5.5)$ и дополнительному
условию $(5.6)$ при условии, что траектории $z$ движения частиц
определены регуляризованным полем скоростей $S_\delta v$ из
уравнения (5.2).

В \cite{dmizv} получена теорема существования решений задачи
$(5.1)-(5.6)$. Как и в случае классической системы Навье-Стокса,
эта теорема утверждает о существовании нелокального сильного
решения для $n=2$ и локального для $n=3$.

{\bf Теорема 5.1}\cite{dmizv}. Пусть $v^0 \in V, \sigma_0 \in
W^1_2(\Omega)^{n^2}, \varphi \in L_2(0,T;L_2(\Omega)^n)$ и
$n=2,3.$ Тогда существует $T_0 > 0$ такое, что для $0<T<T_0$
задача $(5.1)-(5.6)$ имеет хотя бы одно сильное решение$(v
, p)$. При этом $T_0 = \infty$ в случае $n=2$.

Однако, для моделей движения вязкоупругих сред наибольший интерес
представляют не локальные теоремы существования, а глобальные,
поскольку эффект вязкоупругости проявляется, как правило, не
мгновенно, а с течением времени. Ниже будут приведены два
результата о существовании глобальных решений начально-краевых
задач для  регуляризованной и  для исходной модели Джеффриса.

Итак, пусть по--прежнему $\Omega$ - ограниченная область в
$\mathbb{R}^n$ с границей $\partial\Omega$ класса $\mathcal{C}^2,
Q_T =[0,T]\times\Omega$. Далее, кроме стандартных обозначений
$L_q(\Omega), L_q(Q_T), W_q^s(\Omega)$ c  $1\leq q < \infty$ и
$s>0,$ для пространств Лебега и Соболева, состоящих из функций на
$\Omega$ и $Q_T$ со значениями в $\mathbb{R}$, используются
обозначения $L_q(\Omega)^n, L_q(Q_T)^n, W_q^s(\Omega)^n,
B_q^s(\Omega)^n$ для функции со значениями в $\mathbb{R}^n.$ Через
$W_{q,s}^1(\Omega)^{n^2}$ обозначим пространство
$$W_{q,s}^1(\Omega)^{n^2} = \{\sigma = (\sigma_{ij});\sigma \in
W_q^1(\Omega)^{n^2}, \sigma_{ij} =\sigma_{ji},~~~~~~ i,j =
1,2,\ldots,n \}.$$

Через $V_q$ обозначим замыкание $\mathcal{V}$ в норме пространства
$W_q^1(\Omega)^n.$ Обозначим $$W_q^{1,2}(Q_T)=
L_q(0,T;W_q^2(\Omega)^n)\cap L_q(\Omega;W_q^1(0,T)^n).$$ и пусть
$$V_q^{1,2}(Q_T)= W_q^{1,2}(Q_T)\cap L_q(0,T;V_q),$$
\begin{center}
$ P = \{ p: p \in W_2^{0,1}(Q_T), ~~~\int\limits_{\Omega}^{}
p(t,x)dx = 0$ для п. в. $t \in (0,T) \}.$
\end{center}

Ясно, что $V_2^{1,2}(Q_T)=H_s^{1,2}(Q_T)^n.$

Символом $B_q^s(\Omega)$ будем обозначать пространство Бесова
функций на $\Omega$.

Следующий результат нелокальной разрешимости начально-краевой
задачи для  регуляризованной модели Джеффриса при малых правых
частях и начальных данных доказан в работе \cite{DZ1}.

{\bf Теорема 5.2.} Пусть $q \geq 2$ произвольное целое число.
Тогда существуют константы $R_i> 0, i= 1,2,3,$ такие что для любых
$f \in L_q (Q_T)^n, \sigma_0 \in W_{q,s}^1(\Omega)^{n^2},
v^0 \in (B_q^{2-2/q}(\Omega)^n\cap V_q),$ удовлетворяющих
условию
$$\|f\|_{L_q(Q_T)^n}< R_1~~~~\|\sigma_0
\|_{W_{q,s}^1(\Omega)^{n^2}}< R_2~~~~~
\|v^0\|_{B_q^{2-2/q}(\Omega)^n}\leq R_3,$$ задача
(5.1)-(5.6)имеет решение $(v,p)\in V_q^{1,2}(Q_T)\times P.$

Рассмотрим систему уравнений, описывающую не регуляризованную, а
исходную модель Джеффриса
$$\rho(\partial_tv+ \sum\limits_{j=1}^n v_j
\partial x_j v) (t,x) - \mu_1\mathrm{Div} \int\limits_0^t
e^{-\frac{t-s}{\lambda}}\mathcal{E}(v)(s,z(s;t,x))ds
-2\mu_0\mathrm{Div}\mathcal{E}(v)(t,x)+$$
\begin{equation}+ \mathrm{grad} p (t,x) =
e^{-\frac{t}{\lambda}}\mathrm{Div}\sigma_0(z(0;t,x)) + \rho
\varphi(t,x),~~~~~(t,x)\in (0,T)\times\Omega,\end{equation}
\begin{equation} z(\tau;t,x)= x+ \int\limits_t^{\tau}
v(s,z(s;t,x)) ds, ~~~~~~\tau \in [0,T], ~ (t,x) \in
(0,T)\times \Omega. \end{equation}

Рассмотрим для нее также начально-краевую задачу с краевым
условием $(5.4)$ и начальным условием $(5.5)$.

Понятие сильного решения для начально-краевой задачи (5.7),(5.8),
(5.4)-(5.6) аналогично соответствующему понятию для
регуляризованной модели.

Следующий результат нелокальной разрешимости задачи (5.7),(5.8),
(5.4)-(5.6)~ для малых начальных данных и малых правых частей
получен в \cite{DZ1}.

{\bf Теорема 5.3. } Пусть $q>n.$ Тогда существуют константы
$R_i>0, i=1,2,3,
 $ такие, что для любых $f \in L_q(Q_T)^n, \sigma_0 \in W_{q,s}^1(\Omega)^{n^2}, v^0 \in (B_q^{2-2/q}(\Omega)^n\cap
 V_q),$~удовлетворяющих условию $$\|f\|_{L_q(Q_T)^n}<
R_1~~~~\|\sigma_0 \|_{W_{q,s}^1(\Omega)^{n^2}}< R_2~~~~~
\|v^0\|_{B_q^{2-2/q}(\Omega)^n}\leq R_3,$$ задача
(5.7)-(5.8),(5.4)-(5.6) имеет решение $(v, p)\in
V_q^{1,2}(Q_T)\times P.$

\section{Сильные решения начально-краевой задачи для модели Джеффриса
с производной Олдройда}

Начнем этот пункт с наиболее известного результата об этой задаче -- теоремы о существовании и единственности локального
по времени решения начально-краевой задачи для системы уравнений
движения несжимаемой вязкоупругой жидкости с объективной
производной Олдройда в определяющем соотношении, принадлежащую
Гильопе и Со \cite{gs1}.

Пусть $\Omega$ -- ограниченная область (т.е. открытое множество) в
пространстве $\mathbb{R}^3$. Рассмотрим начально-краевую задачу,
состоящую из уравнения движения несжимаемой жидкости
\begin{gather} \frac{\partial v(t,x)}{\partial t} +
\sum\limits_{i=1}^3 v_i(t,x) \frac{\partial v(t,x)}{\partial x_i}
- Div \sigma(t,x) + \nabla p(t,x) = f,\ (t,x)\in
[0,T]\times \Omega,\end{gather} определяющего соотношения
\begin{gather}\sigma(t,x) +\lambda _1\frac{D_a\sigma(t,x)
}{Dt}=2\eta(\mathcal{E}(t,x)+\lambda
_2\frac{D_a\mathcal{E}(t,x)}{Dt}),\ (t,x)\in [0,T]\times
\Omega,\end{gather} условия несжимаемости жидкости
\begin{equation} div
\,v(t,x)=0,(t,x)\in [0,T]\times \Omega,\end{equation} условия
прилипания жидкости к стенкам сосуда
\begin{equation} v(t,x)=0,(t,x)\in [0,T]\times
\partial\Omega,\end{equation}
и начального условия
\begin{equation} v(0,x)=a(x),\sigma(0,x)=\sigma_0(x), x\in
\Omega.\end{equation}

Давление $p$ может быть определено с точностью до константы.
Поэтому к задаче (6.1)-(6.5) обычно для определенности
добавляется условие

\begin{equation}
\int\limits_{\Omega} p(t,x) dx \equiv0
\end{equation}

Напомним необходимые функциональные пространства.
Символом $L_2$ будем обозначать пространство суммируемых с
квадратом на области $\Omega$ функций со значениями в
$\mathbb{R}$, $\mathbb{R}^3$ или пространстве симметричных матриц
$3\times 3$. Символом $H^s$, где $s$ -- натуральное число,
обозначим соболевское пространство функций (со значениями в
$\mathbb{R}$, $\mathbb{R}^3$ или пространстве симметричных матриц
$3\times 3$), суммируемых с квадратом на области $\Omega$ вместе
со своими производными до порядка $s$ включительно; символом
$H^s_0$ -- замыкание в $H^s$ множества бесконечно гладких функций
с компактным в $\Omega$ носителем; символом $H^{-s}$ --
пространство, сопряженное к $H^s_0$.

Теперь мы можем сформулировать результат \cite{gs1}.

\bf Теорема 6.1. \rm Пусть $\Omega$ -- ограниченная область с $C^3$ -
гладкой границей. Пусть $f\in L_2(0,T;H^1)$, $f_t^{'}\in
L_2(0,T;H^{-1})$, $a \in H^2 \bigcap H^1_0$, $div \, a=0$,
$\sigma_0 \in H^1$, $\sigma_0 - \frac{2 \eta \lambda_2}{\lambda_
1}\mathcal{E}(a) \in H^2$.  Тогда существует $t_0>0$ и тройка
$(v,\sigma, p)$ из класса
$$v \in L_2(0,t_0;H^3)\bigcap C([0,t_0];H^2 \cap H^1_0)$$
$$v^{'}_t \in L_2(0,t_0;H^1_0)\bigcap C([0,t_0];L_2)$$
$$p\in L_2(0,t_0;H^2)$$
$$\sigma \in L_2(0,t_0;H^2)\bigcap C([0,t_0];H^1),$$
$$\sigma - \frac{2 \eta \lambda_2}{\lambda_
1}\mathcal{E}(v) \in C([0,t_0];H^2),$$
которая является решением задачи (6.1)-(6.6). Это решение
единственно в указанном классе.

Обозначим через $C_b([0,\infty);X)$ пространство ограниченных на $[0,\infty)$ непрерывных функций со значениями в банаховом пространстве $X$, а через $C^1_b([0,\infty);X)$ пространство функций, имеющих непрерывную ограниченную производную.
При малых данных задача (6.1)-(6.6) имеет глобальное решение \cite{mol1}:

\bf Теорема 6.2. \rm Пусть $\Omega$ -- ограниченная область с $C^3$ -
гладкой границей. Пусть $f\in L_\infty(0,\infty;H^1)$, $f_t^{'}\in
L_\infty(0,\infty;H^{-1})$, $a \in H^2 \bigcap H^1_0$, $div \, a=0$,
$\sigma_0 \in H^1$, $\sigma_0 - \frac{2 \eta \lambda_2}{\lambda_
1}\mathcal{E}(a) \in H^2$, и их нормы малы в указанных пространствах.  Тогда существует и единственна тройка
$(v,\sigma, p)$ из класса
$$v \in L_{2,loc}(0,\infty;H^3)\bigcap C_b([0,\infty);H^2 \cap H^1_0)$$
$$v^{'}_t \in L_{2,loc}(0,\infty;H^1_0)\bigcap C_b([0,\infty);L_2)$$
$$p\in L_{2,loc}(0,\infty;H^2)$$
$$\sigma \in C_b([0,\infty);H^1),$$
$$\sigma - \frac{2 \eta \lambda_2}{\lambda_
1}\mathcal{E}(v) \in C_b([0,\infty);H^2)\bigcap C^1_b([0,\infty);H^1),$$
которая является решением задачи (6.1)-(6.6)  с $T = \infty$.

\section{Существование и единственность сильных решений начальной задачи для системы уравнений движения
нелинейной вязкоупругой среды}

В этом пункте обсуждается задача для уравнений движения для широкого
класса нелинейных вязкоупругих сред. Этот класс содержит, в
частности, ньютоновские и идеальные жидкости, модели Олдройда,
Ларсона, Гизекуса, Фан-Тиена-Таннера, Сприггса, Прандтля, Эйринга
и их комбинации (смеси) \cite{bookour}.

Рассматривается начальная задача для уравнений
движения этого класса нелинейных вязкоупругих сред во всем
пространстве $\mathbb{R}^n, n=2,3$:

\begin{gather}\frac{\partial u}{\partial t}+\sum\limits_{i=1}^{n}
u_i\frac{\partial u}{\partial x_i}=Div\
\mathrm{T}_{\mathrm{H}}+f_0,\ (t,x)\in[0,T]\times\mathbb{R}^n
\\div\ u=0,\ (t,x)\in[0,T]\times\mathbb{R}^n
\end{gather}
\begin{equation} \mathrm{T}_{\mathrm{H}} =
\sigma^s+\sigma^p
\end{equation}
\begin{equation}
\sigma^s=-pI+\Psi(\mathcal{E} (u))
\end{equation}
\begin{gather}
\sigma^p=\sum\limits_{k=1}^r \tau^k\\
\tau^k+\lambda_k\frac{D_{\mathrm{a}}\tau^k}{Dt}+\beta_k(\tau^k,\mathcal{E})=2
\eta_k \mathcal{E}
\end{gather}

Здесь $u$ --- вектор скорости,  $f_0$ --- поле внешних сил,
$\mathrm{T}_{\mathrm{H}}$ --- тензор напряжений (все они зависят
от точки $x$ пространства $\mathbb{R}^n$, $n=2,3$ и момента
времени $t$), $r$ --- натуральное число, $\tau^k, \sigma^s$,
$\sigma^p$ --- составляющие тензора напряжений,
$\mathcal{E}(u)=\frac {1}{2}(\nabla u+\nabla u^T)$ -- тензор
скоростей деформации, $p$ --- давление, $\lambda_k>0$ -- времена
релаксации, $\eta_k>0$ -- вязкости, $k=1,\dots,r$. Выражение
$\frac{D_{\mathrm{a}} A}{Dt}$ есть объективная
(олдройдовская) производная тензора.
Кроме того, $\Psi$ и $\beta_k$ в (7.4) и (7.6) это
известные нелинейные функции со значениями во множестве матриц
$n\times n$, которые имеют следующий вид:
\begin{gather}
\Psi(\mathcal{E})= \varphi_1\mathcal{E} +
\varphi_2\mathcal{E}^2\\
\beta_k(\tau, \mathcal{E})=\alpha^k_0 I+\alpha^k_1 \mathcal{E}+
\alpha^k_2 \mathcal{E}^2 + \alpha^k_3 \tau+ \alpha^k_4 \tau^2
+\alpha^k_5 (\mathcal{E}\tau+\tau\mathcal{E})+\notag\\+ \alpha^k_6
(\mathcal{E}^2\tau+\tau\mathcal{E}^2)+\alpha^k_7
(\mathcal{E}\tau^2+\tau^2\mathcal{E})+\alpha^k_8
(\mathcal{E}^2\tau^2+\tau^2\mathcal{E}^2)
\end{gather} где $\varphi_1$, $\varphi_2$ и $\alpha^k_j$ суть произвольные скалярные
функции следующих аргументов
$$\varphi_i=\varphi_i(Tr(\mathcal{E}^2), det \mathcal{E}), i=1,2$$
$$\alpha^k_j=\alpha^k_j(Tr \mathcal{E}^2,Tr
\mathcal{E}^3,Tr(\tau),
Tr(\tau^2),Tr(\tau^3),Tr(\tau\mathcal{E}),$$ $$
Tr(\tau^2\mathcal{E}),Tr(\tau\mathcal{E}^2),Tr(\tau^2\mathcal{E}^2)),\
k=1,\dots,r; j=0,\dots,8$$

Давление $p$ вообще может быть определено с точностью до
константы. Для определенности накладывается условие

\begin{equation}
\int\limits _ {\Omega} p (t,x) dx \equiv0
\end{equation}
где $\Omega$ - фиксированная ограниченная область в
$\mathbb{R}^n$.

Начальные условия имеют вид:
\begin{equation}
\begin{array}{l}
u(0,x)=a(x),\ \tau^k(0,x)=\tau^k_0(x),\ x \in \mathbb{R}^n,
k=1,\dots,r
\end{array}
\end{equation}

Обозначим $\eta_0=\frac {\varphi_1(0,0)} 2$. Будем предполагать,
что $\eta_0>0$, а $\varphi_1$, $\varphi_2$, $\alpha^k_j$
соответственно $C^4$-,  $C^4$- и $C^3$-гладкие функции и
$$\alpha^k_0(\theta)=\alpha^k_1(\theta)=\alpha^k_3(\theta)=0, \frac {\partial \alpha^k_0(\theta)}{\partial{Tr(\tau)}}=0$$ ($\theta$ обозначает
точку $(0,0,0,0,0,0,0,0,0)$).

 Мы будем использовать
функциональные пространства типа Соболева \linebreak $H_V^m=\{u
\in H^m(\mathbb{R}^n,\mathbb{R}^n), \ div\, u=0\}$ и
$H_M^m=H^m(\mathbb{R}^n, \mathbb{R}_S^{n \times n})$.

\bf Теорема 7.1. \rm \cite{vz3}.  Пусть $a\in H_V^3, \ \tau^k_0 \in H_M^3,$
$k=1, \dots, r$, $f_0 \in L_1 (0, T; H^3 (\mathbb {R} ^n,
\mathbb {R} ^n))\cap L_2 (0, T; H^2 (\mathbb {R} ^n, \mathbb {R}
^n))$. Тогда существует такая константа $K_0>0$, не зависящая от
$T$, что при $$\|a\|_3+\sum\limits_{k=1}^r
\|\tau_0^k\|_3+\|f_0\|_{L_1(0,T;H^3)}<K_0$$ задача (7.1)-(7.10) имеет
решение в классе
$$u\in L_2(0,T;H_V^4) \cap C([0,T];H_V^3)\cap W_2^1(0,T;H_V^2)$$
$$
\mathrm{T}_{\mathrm{H}} \in L_2(0,T;H_{M,loc}^{3})$$ $$ p\in
L_2(0,T;H^3_{loc}(\mathbb{R}^n,\mathbb{R}))$$ Кроме того, о
составляющих тензора напряжений имеется следующая информация:
$$\sigma^s+pI\in L_2(0,T;H_M^3) \cap C([0,T];H_M^2)\cap
W_2^1(0,T;H_M^1)$$ $$ \sigma^p, \tau^k \in L_\infty(0,T;H_M^3)\cap
C([0,T];H_M^2)\cap C^1([0,T];H_M^1)$$ $$ k=1,\dots,r $$ Это
решение единственно в этом классе.

Если $f_0 \in L_1 (0, + \infty; H^3 (\mathbb {R} ^n, \mathbb {R}
^n)) \cap L_2 (0, + \infty; H^2 (\mathbb {R} ^n, \mathbb {R} ^n))$
и $$\|a\|_3+\sum\limits_{k=1}^r
\|\tau_0^k\|_3+\|f_0\|_{L_1(0,+\infty;H^3)}<K_0$$ то задача
(7.1)-(7.10) имеет единственное решение в указанном классе при всяком
$T>0$. \rm

\section{О нерешенных проблемах} Основные нерешенные проблемы по задаче (0.2),(0.4),(0.5),(2.18) можно условно разбить на две группы: a) проблемы, имеющие аналогичную нерешенную проблему для системы Навье-Стокса\footnote{Хотя потенциально решение такой проблемы может быть основано на свойствах "вязкоупругой" проблемы, не связано с соответствующими проблемами для Навье-Стокса и метод ее решения не обязательно подойдет для Навье-Стокса.}, и б) проблемы, характерные сугубо\footnote{Т.е., как правило, аналогичная проблема для Навье-Стокса решена.} для (0.2),(0.4),(0.5),(2.18).  

К первой группе относятся проблема сильной глобальной разрешимости начально-краевой задачи, проблема регулярности слабых решений,  проблема единственности слабых решений, оценка хаусдорфовой размерности траекторных и глобальных аттракторов\footnote{Наличие которых пока не установлено.} (все эти задачи требуют решения без ограничений на коэффициенты, внешнюю силу и начальные данные).  
Важно, что для Навье-Стокса эти задачи актуальны при $n \geq 3$, а для (0.2),(0.4),(0.5),(2.18) --- при $n \geq 2$. 

В рамках второй группы важен блок задач о построении хоть каких-нибудь обобщенных решений стационарной краевой задачи и эволюционной\footnote{Включая вопрос (см. п.3) о проверке справедливости результатов \cite{lmas} в коротационном случае (3.2).} начально-краевой задачи, а также задача о наличии периодических по времени решений,  все --- без ограничений на коэффициенты, внешнюю силу и начальные данные. Для эволюционной задачи речь идет о глобальных по времени решениях. В случае успешного (даже частичного) решения задачи о существовании глобальных решений встает вопрос о построении любых, пусть обобщенных, аттракторов и/или \emph{инерционных многообразий} (о последнем понятии см. \cite{tem2}).

\end{document}